\documentclass[letterpaper, 10 pt, conference]{ieeeconf}  

\IEEEoverridecommandlockouts                              
\usepackage{graphicx} 
\usepackage{amsmath} 
\usepackage{amssymb}  

\usepackage{color}
\usepackage{hyperref}
\usepackage{cite}

\title{\LARGE \bf
Singularly perturbed phase response curves
}

\author{Pierre Sacr\'e and Alessio Franci
\thanks{P. Sacr\'e is with the Department of Biomedical Engineering, 
		Johns Hopkins University, Baltimore, MD
        {(\href{mailto:p.sacre@jhu.edu}{p.sacre@jhu.edu})}.}%
\thanks{A. Franci is with the Department of Electrical Engineering and Computer Science, University of Li\`{e}ge, Belgium, and with the Department of Engineering, Cambridge University, United Kingdom  
        {(\href{mailto:af529@cam.ac.uk}{af529@cam.ac.uk})}.}%
}

\newtheorem{remark}{Remark} 

\newcommand{\field}[1]{\mathbb{#1}} 
 
\newcommand{\Rbb}{\field{R}}

\newcommand{\Sbb}{\field{S}}

\newcommand{\eqdef}{\mathrel{:=}}
\newcommand{\reveqdef}{\mathrel{=:}}

\makeatletter
\newcommand{\prap}[1]{\allowbreak\if@display\mkern18mu\else\mkern8mu\fi(#1)}
\newcommand{\pwrap}[1]{\prap{{\operator@font wrapped\;in}\mkern6mu#1}}
\makeatother

\newcommand{\zeroinput}{\mathbf{0}}

\newcommand{\fvecf}{f}
\newcommand{\fvecs}{g}

\newcommand{\flow}{\phi}
\newcommand{\perorb}{\gamma}
\newcommand{\persol}{x^{\perorb}}

\newcommand{\PRC}{Q}

\newcommand{\iPRCu}{q}

\newcommand{\ech}{1_{+}}



\newcommand{\critman}{\mathcal{S}}

\newcommand{\eqman}{\mathcal{S}}

\newcommand{\stableeqm }{\eqman^{a}_{-}}
\newcommand{\stableeqp }{\eqman^{a}_{+}}
\newcommand{\stableeqb }{\eqman^{a}_{\bullet}}

\newcommand{\unstableeq}{\eqman^{r}}

\newcommand{\foldp}{\mathcal{F}_{+}}
\newcommand{\foldm}{\mathcal{F}_{-}}

\newcommand{\critfun}{b}

\newcommand{\perfast}{T_{\text{f}}}
\newcommand{\perslow}{T_{\text{s}}}

\newcommand{\freqfast}{\omega_{\text{f}}}
\newcommand{\freqslow}{\omega_{\text{s}}}

\newcommand{\Deltafast}{\Delta_{\text{f}}}
\newcommand{\Deltaslow}{\Delta_{\text{s}}}

\graphicspath{{./_figures/pdf/}}

\newcommand{\myfigurename}{Figure}

\newcommand{\flowf}{\flow_{\text{f}}}

\newcounter{condition}
\begin{document}

\fontencoding{T1}\fontsize{10}{12pt}\selectfont

\maketitle
\thispagestyle{empty}
\pagestyle{empty}

\begin{abstract} 
In this paper we propose a novel geometric method, based on singular perturbations, to approximate isochrones of relaxation oscillators and predict the qualitative shape of their (finite) phase response curve. This approach complements the infinitesimal phase response curve approach to relaxation oscillators and overcomes its limitations near the singular limit. We illustrate the power of the methodology by deriving semi-analytic formula for the (finite) phase response curve of generic planar relaxation oscillators to impulses and square pulses of finite duration and verify its goodness numerically on the FitzHugh--Nagumo model.
\end{abstract}

\section{Introduction}

The phase response curve (PRC) characterizes the input--output behavior of oscillatory systems \cite{Winfree:1980ue,Glass:1988ub}. It has wide applications ranging from oscillator control \cite{Efimov:2009fr,Danzl:2009co} to the analysis of oscillator network synchronization \cite{Mauroy:2012vi,Dorfler:2013fk}. Systematic and analytic prediction of an oscillator phase response curve is a hard task in general and it can be accomplished only in very specific cases. This usually leads to intense numerical investigations, which might weaken the relevance of phase response curve approach in engineering applications.

The classical approach relies on analytically computing the infinitesimal (linearized) phase response curve and then use convolution to compute the phase response curve for generic inputs~\cite{Sacre:2014a a}. However, because the measure of the region where the linear approximation is valid shrinks to zero as the time-scale separation is increased, this approach is valid only for inputs of infinitesimally small amplitude~\cite{Izhikevich:2000hb}.

To overcome this limitations, we used a complementary approach and study the global structure of the oscillator isochrones from a geometric viewpoint, using singular perturbation method \cite{Krupa:2001ez}. Based on this analysis we derived semi-analytic formulas to approximate the (finite) phase response curve to arbitrary inputs in the singular limit. 
As opposed to the infinitesimal phase response curve approach, the error between the real and the predicted singular phase response curve goes to zero as the time-scale separation increases. 
The method is thoroughly illustrated  by predicting the phase response curve of a generic relaxation oscillator to impulses and square pulses of finite duration. 
These results are a first step toward a geometric theory for (finite) phase response curves of singularly perturbed oscillators, including  bursters~\cite{Franci:2013uq}.

The results of the article are primarily drawn from the Ph.D. dissertation of the first author \cite{Sacre:2013ys}.

This paper is organized as follows.
Section~\ref{sec:relax_oscill} summarizes basics concepts of singular perturbation theory and describes the underlying geometry  of relaxation oscillator dynamics.
Section~\ref{sec:limitations} stresses the limitation of the  \emph{infinitesimal} phase response curve approach in the context of relaxation oscillators with inputs acting on the fast variable dynamics.
Section~\ref{sec:singular_prc} introduces the novel concept of singularly perturbed phase response curves predicted from the singular limit.
Section~\ref{sec:application} illustrates our geometric approach on the FitzHugh--Nagumo neural model.%


\section{Relaxation oscillators and their geometry} \label{sec:relax_oscill}

We consider a  two-dimensional fast-slow dynamical system of the form 
\begin{subequations} \label{eq:relax_fast_syst}
\begin{align}
	\dot{x} & = \phantom{\epsilon\,}\fvecf(x) - z + u,\label{eq:relax_fast_syst_f} \\
	\dot{z} & = \epsilon\,\fvecs(x,z),  \label{eq:relax_fast_syst_g}
\end{align}
\end{subequations}
where $\dot{}$ denotes differentiation with respect to the time~$t$, $(x,z)\in\Rbb^2$, $u\in\Rbb$, and $0<\epsilon \ll 1$.
The solution at time~$t$ to the initial value problem \eqref{eq:relax_fast_syst} from the initial condition $(x_0,z_0)\in\Rbb^2$ at time~$0$ is denoted by $\flowf^\epsilon(t,(x_0,z_0),u(\cdot))$, with $\flowf^\epsilon(0,(x_0,z_0),u(\cdot)) = (x_0,z_0)$.
In the slow time scale $\tau\eqdef\epsilon\,t$, dynamics \eqref{eq:relax_fast_syst} become
\begin{subequations} \label{eq:relax_slow_syst}
\begin{align}
	\epsilon\,x' & = \fvecf(x) - z + u, \label{eq:relax_slow_syst_x}\\
	z' & =\fvecs(x,z),\label{eq:relax_slow_syst_z}
\end{align}
\end{subequations}
where $'$ denotes differentiation with respect to the slow time~$\tau$.  For $\epsilon\neq0$, \eqref{eq:relax_fast_syst} and \eqref{eq:relax_slow_syst} are equivalent. We call \eqref{eq:relax_fast_syst} the fast dynamics and \eqref{eq:relax_slow_syst} the slow dynamics.
In the limit $\epsilon\to 0$, we obtain from \eqref{eq:relax_fast_syst} and \eqref{eq:relax_slow_syst} the \emph{layer} dynamics
\begin{subequations} \label{eq:relax_fast_syst layer}
\begin{align}
	\dot{x} & = \fvecf(x) - z + u,\label{eq:relax_fast_syst_f layer} \\
	\dot{z} & =0, \label{eq:relax_fast_syst_g layer}
\end{align}
\end{subequations}
describing the fast evolution far from the critical manifold $\critman^0\eqdef\left\{(x,z)\in\Rbb^2:\fvecf(x) - z + u=0\right\}$, and the \emph{reduced} dynamics
\begin{subequations} \label{eq:relax_slow_syst reduced}
\begin{align}
	0& = \fvecf(x) - z + u, \label{eq:relax_slow_syst_x reduced}\\
	z' & =\fvecs(x,z),\label{eq:relax_slow_syst_z reduced}
\end{align}
\end{subequations}
describing the slow evolution along $\critman^0$. 

\begin{figure}
	\centering
	\includegraphics[scale=0.8]{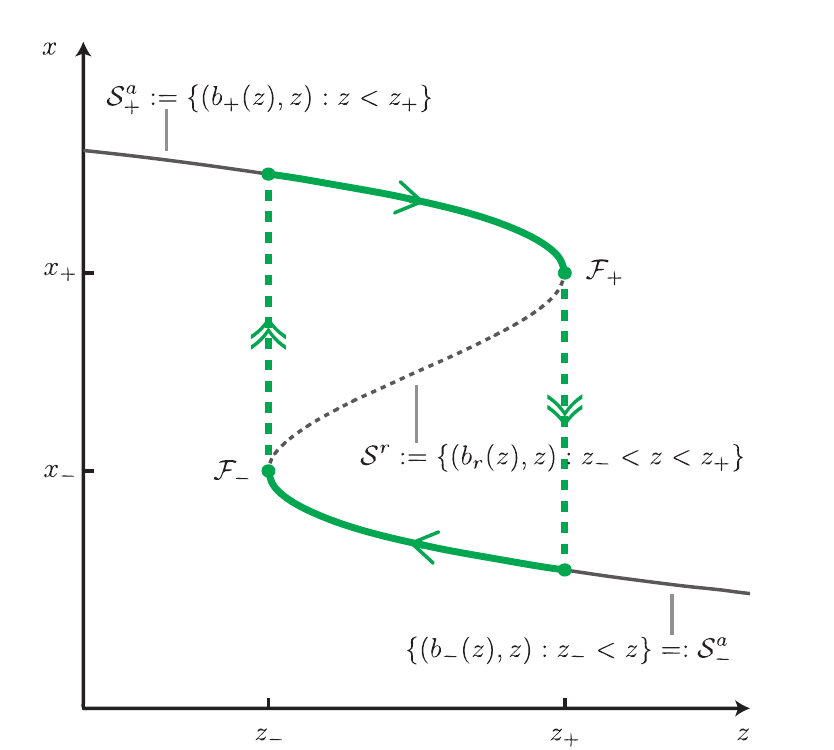}
	\caption{%
		Geometry of relaxation oscillators. 
		The critical manifold $\critman^0$ is a \mbox{S-shaped} curve.
		Under some technical assumptions \cite{Krupa:2001ez}, the singular system~\eqref{eq:relax_fast_syst} admits a singular periodic orbit~$\perorb^0$ defined as the union of two pieces of the critical manifold associated with a slow evolution (green solid lines) and two critical fibers associated with jumps (green dashed lines).
		}
	\label{fig:relax_oscillator_geom_merge}
\end{figure}

Under some mild technical assumptions \cite[Theorem 2.1]{Krupa:2001ez}, in particular that the critical manifold $\critman^0$ is S-shaped, the zero-input system~\eqref{eq:relax_fast_syst} has a unique periodic orbit $\perorb^\epsilon$ sliding along the stable branches of $\critman^0$ and shadowing the singular periodic orbit $\perorb^0$ illustrated in~\myfigurename~\ref{fig:relax_oscillator_geom_merge}. The singular periodic orbit~$\perorb^0$ is defined as the union of two pieces of the critical manifold associated with a slow evolution and two critical fibers associated with jumps.

\begin{remark}
In the slow time scale, the singularly perturbed period~$\perslow^\epsilon$ converges towards the singular period~$\perslow^0$, which is equal to the finite time required to slide along both portions of the critical manifold (jumps are instantaneous), that is, $\lim_{\epsilon \rightarrow 0} \perslow^\epsilon \reveqdef \perslow^0$.
	In the fast time scale, the singularly perturbed period~$\perfast^\epsilon$ blows up to infinity, that is, $\lim_{\epsilon \rightarrow 0} \perfast^\epsilon \reveqdef\perfast^0$, with $\lim_{\epsilon \rightarrow 0} \perfast^\epsilon = \lim_{\epsilon \rightarrow 0} \perslow^\epsilon/\epsilon= +\infty$.
	Note that corresponding angular frequencies are denoted $\freqfast^\epsilon = 2\pi/\perfast^\epsilon$ and $\freqslow^\epsilon = 2\pi/\perslow^\epsilon$.
\end{remark}

\section{Limitation of the infinitesimal approximation \\for finite phase response curves} \label{sec:limitations}

In this section, we  introduce the concepts of phase map and phase response curves following the terminology of~\cite{Winfree:1980ue} and~\cite{Glass:1988ub}. 
(The interested reader is referred to \cite{Sacre:2014aa} for detailed.)
We~emphasize the limitation of the \emph{infinitesimal} approximation for   \emph{finite} phase response curves of relaxation oscillators.

\subsection{Phase map and isochrons}
Because of the periodic nature of its steady-state behavior, it is appealing to study the oscillator dynamics  on the unit circle $\Sbb^1$. The key ingredients of this phase reduction are the concepts of phase map and isochrons.
 
The (asymptotic) \emph{phase map} $\Theta^\epsilon:\mathcal{B}(\gamma^\epsilon)\subseteq\Rbb^2\rightarrow\Sbb^1$ is a mapping that associates with every point in the basin of attraction $\mathcal{B}(\gamma^\epsilon)$ a phase on the unit circle $\Sbb^1$. 
It is defined such that the phase variable $\theta(t)\eqdef\Theta^\epsilon(\flowf^\epsilon(t,(x_0,z_0),u(\cdot)))$, that is, the image of the flow through the phase map, linearly increases  with time for the input~$\zeroinput$. 

\emph{Isochron} $\mathcal{I}^\epsilon(\theta)$ are the set of points in $\mathcal{B}(\gamma^\epsilon)$ that are associated with the same asymptotic phase $\theta$ on the unit circle~$\Sbb^1$. Isochrons are level sets of the phase map. 

\subsection{Phase response curves}

An input $u(\cdot)$ is  \emph{phase-resetting} if the  solution of~\eqref{eq:relax_fast_syst} forced by $u(\cdot)$ asymptotically converges to the periodic orbit. 

The \emph{(finite) phase response curve} $Q^\epsilon(\cdot;u(\cdot)):\mathbb S^1\to [-\pi,\pi)$ associates with each phase the asymptotic phase shift of system \eqref{eq:relax_fast_syst} in response to a phase-resetting input $u(\cdot)$.

The \emph{infinitesimal phase response curve} $q^\epsilon:\mathbb S^1\to\mathbb R$ is the relative asymptotic phase shift of system \eqref{eq:relax_fast_syst} in response to an infinitesimal phase-resetting impulse (Dirac $\delta$ function), that is, $q^\epsilon(\cdot) \eqdef \lim_{\alpha\rightarrow0} Q^\epsilon(\cdot;\alpha\,\delta(\cdot))/\alpha$.

\subsection{Limitation of the infinitesimal approximation}

In the classical approach, the (finite) phase response curve is approximated by the ``convolution'' between the infinitesimal phase response curve and the phase-resetting input 
\begin{equation*}
	\PRC^{\epsilon}(\cdot;u(\cdot)) =  \underbrace{\int_{0}^{+\infty}\iPRCu^{\epsilon}(\omega\,s+\cdot) \, u(s) \,ds}_{\eqdef \PRC^{\epsilon}_{\text{inf}}(\cdot;u(\cdot))} + \mathcal{O}(\|u(\cdot)\|^2).
\end{equation*}
For relaxation oscillators, this approximation is only valid for inputs that are much smaller than the singular perturbation parameter, that is, $0 < \|u(\cdot)\| \ll \epsilon \ll 1$ (see \cite{Izhikevich:2000hb} for details).
Therefore, the domain of validity of this approximation vanishes in the singular limit ($\epsilon \rightarrow 0$). 

Intuitively, this limitation comes from the fact that, the singular trajectory induced by the input $u(\cdot)$ might jump instantaneously from one branch of the critical manifold to the other. This behavior involves a global phenomenon that cannot be captured by a local approximation.

\section{Prediction from the singular limit \\for finite phase response curves} \label{sec:singular_prc}

The main idea of our  approach is to take advantage of time-scale separation to study the \emph{finite} phase response curve in the singular limit. For a sufficiently small singular parameter $0<\epsilon\ll1$, the \emph{singularly perturbed} finite phase response curve~$\PRC^{\epsilon}(\cdot;u(\cdot))$ can be approximated by a \emph{singular} finite phase response curve~$\PRC^{0}(\cdot;u(\cdot))$, that is,
\begin{equation*}
	\PRC^{\epsilon}(\cdot;u(\cdot)) = \PRC^{0}(\cdot;u(\cdot)) + \mathcal{O}\left( \epsilon^{\beta} \right),
\end{equation*}
for any phase-resetting input $u(\cdot)$ and with $0 < \beta \leq 1$. Geometric singular perturbation arguments let $\beta\sim1/2$ \cite{Krupa:2001ez}. 

\begin{figure}
	\centering
	\includegraphics[scale=0.8]{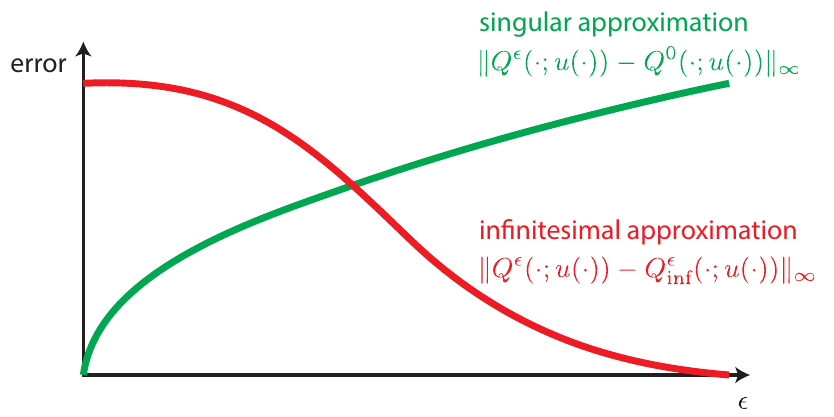}
	\caption{Qualitative trade-off between the infinitesimal approximation and the singular approximation as a function of the time-scale separation $\epsilon$.}
	\label{fig:trade-off}
\end{figure}

\myfigurename~\ref{fig:trade-off} illustrates qualitatively the existing trade-off between the infinitesimal approximation and the singular approximation as a function of the time-scale separation $\epsilon$.

In this section, we show how to exploit the geometry of relaxation oscillators to approximate the phase map near the singular limit and use this construction to predict the shape of the (finite) phase response curve where, for the sake of illustration, we stick to impulses, that is, $u(\cdot) = \alpha\,\delta(\cdot)$, and square pulses of finite duration, that is, $u(\cdot) = \bar{u}\,\left[\ech(\cdot) - \ech(\cdot - \Delta)\right]$.

\subsection{Singular phase map and isochrons}

\begin{figure}
	\centering
	\includegraphics[scale=0.8]{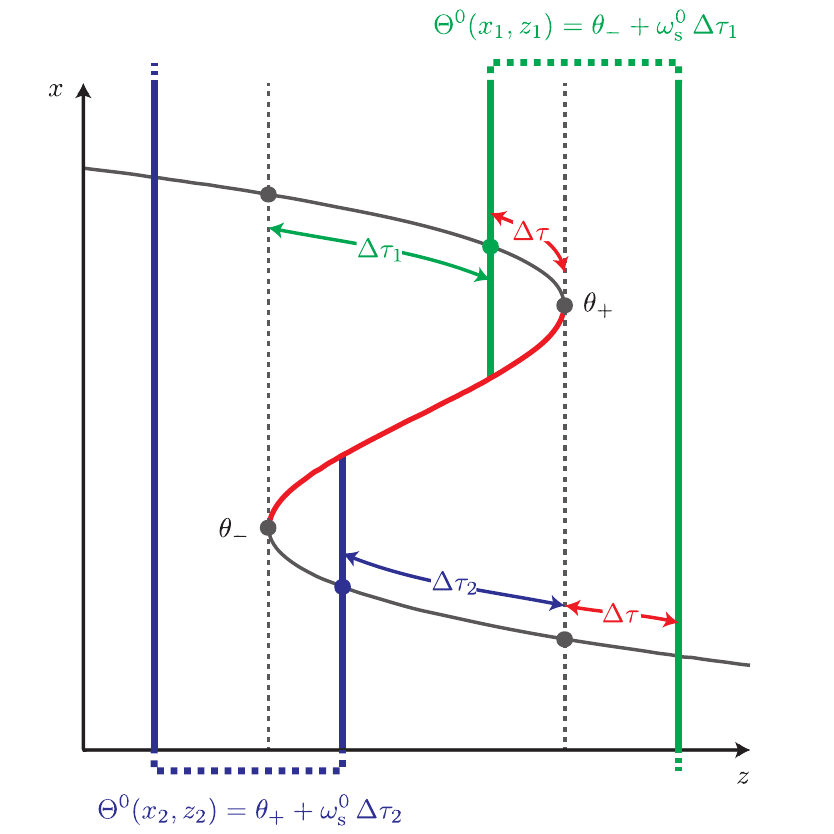}
	\caption{%
		Geometric construction of (asymptotic) singular phase map.
		The phase map associates with each point on the periodic orbit a phase which corresponds to the normalized time $\freqslow^0\,\Delta\tau$ required to reach this point from the reference position $(x_-,z_-)$. For points  on the lower branch, it is convenient to measure the  normalized time from $(x_+,z_+)$ and to add the phase $\theta_+ \eqdef \freqslow^0\,\Delta\tau_{+}$.
		Because all points on a same vertical ray (in the bistable region) and converging to the same branch instantaneously jump on the branch in the singular limit, the asymptotic phase map associates them with the same asymptotic phase. In addition, other vertical lines (outside the bistable region) are associated with the same phase because these points converge in the same $\Delta\tau \pmod{\perslow^0}$ to $(x_+,z_+)$.
		}
	\label{fig:relax_oscillator_iso}
\end{figure}

A first step towards the prediction of singular (finite) phase response curves is the construction of the (asymptotic) phase map and  isochrons for the system \eqref{eq:relax_fast_syst} in the singular limit.
Their construction relies on the fast-slow dynamics \eqref{eq:relax_fast_syst layer}--\eqref{eq:relax_slow_syst reduced} and is illustrated in \myfigurename~\ref{fig:relax_oscillator_iso}.

Since the singular periodic orbit~$\perorb^0$ is a one-dimensional piece-wise smooth curve in $\Rbb^2$, it is naturally parameterized in terms of a single scalar phase on the unit circle~$\Sbb^1$. As in the nonsingular case, the phase map will be chosen such that the phase variable linearly increases with time. 

We choose to associate the zero-phase reference position on the singular periodic orbit with the lower fold $(x_{-},z_{-})$, that is $\Theta^0(x_{-},z_{-})\reveqdef \theta_{-} = 0$. As jumps are instantaneous in the singular limit, all points of the (weakly) unstable critical fiber joining $(x_{-},z_{-})$ to $(\critfun_{+}(z_-),z_-)$ are also associated with a phase equal to zero.

The phase $\theta$ associated with a point~$(x,z)$ is the ``normalized'' fraction of (slow) time $\freqslow^0\,\Delta \tau$  needed to reach this point along the reduced dynamics \eqref{eq:relax_slow_syst reduced} flow from the reference initial condition.
For a point $(x_1,z_1)$ on the upper branch, the phase will be given by 
\begin{equation*}
	\Theta^0(x_1,z_1) \eqdef  \freqslow^0\,\Delta\tau_1.
\end{equation*}
For a point $(x_2,z_2)$ on the lower branch, the phase will be given by
\begin{equation*}
	\Theta^0(x_2,z_2) \eqdef 
	\freqslow^0\,\Delta\tau_{+} + \freqslow^0\,\Delta\tau_{2} 
\end{equation*}
where the first term corresponds to the flowing time on the upper branch (up to the upper fold) and the second term corresponds to the remaining flowing time on the lower branch.
To simplify notation, it is convenient to denote by $\Theta^0(x_{+},z_{+}) \reveqdef \theta_{+} = \freqslow^0\,\Delta\tau_{+}$ the phase associated with the upper fold (and all points of the (weakly) unstable critical fiber joining $(x_{+},z_{+})$ to $(\critfun_{-}(z_+),z_+)$).

The notion of singular phase map can be extended to any point $(x,z)$ in the basin of attraction of the singular periodic orbit. 
Because, in the singular limit, any singular trajectory starting from $(x,z)$ instantaneously jumps from its initial condition to a branch of the critical manifold, all points on the same vertical line (that is, with the same value of slow variable $z$) jumping to the same branch are associated with the same phase.

\begin{itemize}
	\item All points on the line $z = z_{-}$ (resp. $z = z_{+}$) are associated with the phase $\theta_{-}$ (resp. $\theta_{+}$).
	\item For points with a slow variable in the bistable range, the asymptotic phase $\theta_1$ of a point $(x_1,z_1)$ belonging to the basin of attraction of the upper (resp.~lower) branch is thus given by the phase $\theta_1$ of the point at the intersection between the line $z=z_1$ and the upper (resp. lower) branch of the singular periodic orbit $\perorb^0$.
	\item In addition, all points outside the bistable range that converge to  the upper fold in the same time interval $\Delta\tau \pmod{\perslow^0}$ as $(x_1,z_1)$ are also associated with the asymptotic phase $\theta_1$.
\end{itemize}

An elegant way to summarize the definition of the singular (asymptotic) phase map is 
\begin{equation*}
	\Theta^0(x,z) =
	\begin{cases}
			\theta_{-} + \freqslow^0 \, \psi_{+}(z_-,z,0) \pmod{2\pi}, & \text{if (C1)}, \\
			\theta_{+} + \freqslow^0 \, \psi_{-}(z_+,z,0) \pmod{2\pi}, & \text{if (C2)},
	\end{cases}
\end{equation*}
with 
\begin{align*}
	\text{(C1)} & \equiv (x,z)\in\mathcal{B}(\stableeqp) \cup \foldp, \\
	\text{(C2)} & \equiv (x,z)\in\mathcal{B}(\stableeqm) \cup \foldm.
\end{align*}
where $\psi_{\bullet}(z_0,z_\tau,\zeroinput)$ (with $\bullet$ standing for $+$ or $-$) are functions that measure the time needed to travel along the critical manifold from the initial condition $z_0$ to finial condition $z_{\tau}$ and
where $\mathcal{B}(\stableeqb)$ is the set of points that jumps to the stable branch $\stableeqb$ of the critical manifold.

Singular isochrons are thus vertical lines for values of~$z$ outside the bistable range and vertical rays for values of~$z$ inside the bistable range. In the bistable region, vertical rays are separated by the repulsive branch $\unstableeq$ of the critical manifold. The vertical ray and the vertical lines associated with the same phase join at infinity (see \myfigurename~\ref{fig:relax_oscillator_iso}).

For constant inputs $u(\cdot)\equiv\bar{u}$, the function $\psi_{\bullet}(z_0,z_\tau,\bar{u})$ can easily be computed by integrating the reduced dynamics~\eqref{eq:relax_slow_syst reduced}  on the stable branches of the critical manifold and they read
\begin{equation*}
	\psi_{\bullet}(z_0,z_\tau,\bar{u}) = \int_{z_0}^{z_\tau} \frac{1}{g(\critfun_{\bullet}(\xi-\bar{u}),\xi)} d\xi .
\end{equation*}

\begin{remark}
	For presentation convenience, we intentionally do not consider the unstable branch of the critical manifold~$\unstableeq$ as being part of the basin of attraction of the singular periodic orbit. For small $\epsilon$, this repulsive branch is perturbed into a repulsive set which has zero Lebesgue measure.
\end{remark}

\begin{remark}
	For convenience, the singular periodic orbit~$\perorb^0$ is parameterized by the map $\persol:\Sbb^1\rightarrow\perorb^0$ that associates with each phase $\theta$ on the unit circle a point $(\persol(\theta),z^\gamma(\theta))$ on the periodic orbit.
\end{remark}

\subsection{Singular (finite) phase response curves}

We derive the singular (finite) phase response curve for two inputs: impulses, that is, $u(\cdot) = \alpha\,\delta(\cdot)$, and square pulses of finite duration, that is, $u(\cdot) = \bar{u}\,\left[\ech(\cdot) - \ech(\cdot - \Delta)\right]$.

\subsubsection{Impulse}
An impulse $u(\cdot) = \alpha\,\delta(\cdot)$ induces a jump of the fast variable $x$ in the fast-slow dynamics \eqref{eq:relax_fast_syst layer}--\eqref{eq:relax_slow_syst reduced}. The singular (finite) phase response curve is thus given by
\begin{equation*}
	\PRC^0(\theta;\alpha\,\delta(\cdot)) = \Theta^0(x^\gamma(\theta),z^\gamma(\theta) + \alpha) - \theta.
\end{equation*}

\begin{figure}
	\centering
	\includegraphics[scale=0.8]{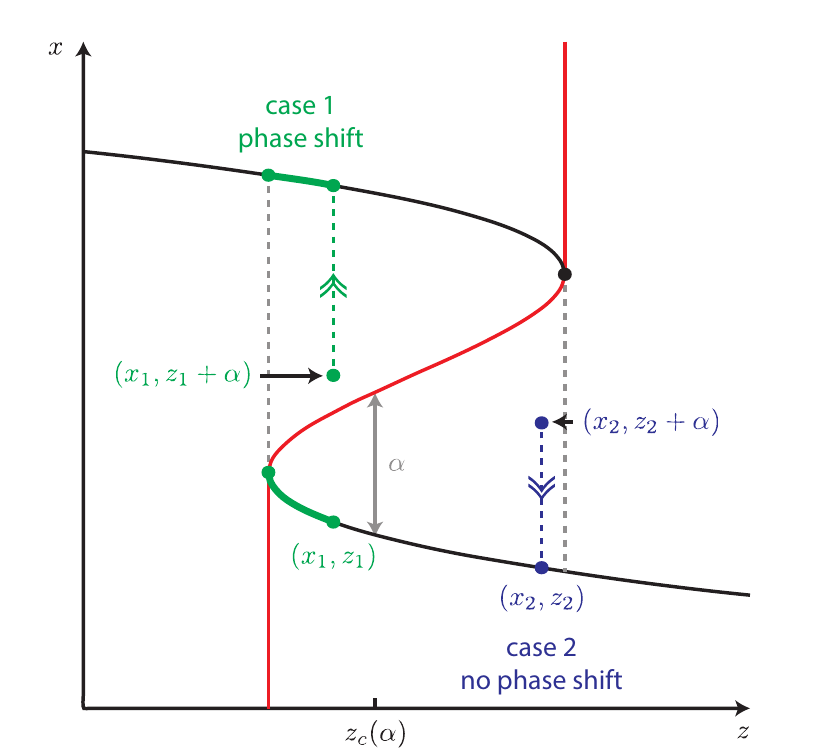}
	\caption{
		Effect of positive impulses in the fast-slow dynamics \eqref{eq:relax_fast_syst layer}--\eqref{eq:relax_slow_syst reduced}.
		(Case 1)~Close enough to the lower fold (on the lower branch), the reset state crosses the separatrix (red curve) and converges toward the upper branch instantaneously. The phase shift corresponds to the phase difference corresponding to the skipped portions of the singular periodic orbit  (green).
		(Case 2)~Far from the lower fold (on the lower branch) or on the upper branch, the reset state converges back to the initial state instantaneously. As a consequence, no phase shift is produced. 
		}
	\label{fig:relax_oscillator_impulse}
\end{figure}

As illustrated on \myfigurename~\ref{fig:relax_oscillator_impulse}, if the impulse lets the state cross the unstable branch of the critical manifold (case~1), it produces a phase shift. In the opposite case (case~2), the state converges back to the initial condition almost instantaneously.

For simplicity, we assume monotonicity of this separatrix in the bistable region (that is, $(\partial \critfun_{r}/\partial z)(z)>0$).

Given a positive impulse of amplitude~$\alpha$, there exists a critical value $z_{c}(\alpha)$ of the slow variable such that a trajectory starting on the lower branch crosses the separatrix under the effect of the impulse for all $z$, such that $z_{-} \leq z < z_{c}(\alpha)$.
The critical value $z_{c}(\alpha)$ is given  by
\begin{equation*}
 	z_{c}(\alpha) = \{ z\in\Rbb :  \critfun_{-}(z) + \alpha = \critfun_{r}(z) \}. 
\end{equation*}
The asymptotic phase associated with this critical point $(\critfun_{-}(z_{c}(\alpha)),z_{c}(\alpha))$ on the stable branch is denoted~$\theta_c(\alpha)$. The phase shift~$\Delta\theta$ induced by an impulse corresponds to the portion of singular periodic orbit skipped due to the impulse.

The phase response curve is given by
\begin{align*}
	\PRC^0(\theta;\alpha\,\delta(\cdot)) &= 
	\begin{cases}
		\theta_{-} + \freqslow^0 \, \psi_{+}(z_{-},z^{\perorb}(\theta),0) - \theta, &\text{if (C3)}, \\
		0, & \text{o/w},
	\end{cases}
\end{align*}
where (C3) stands for $\theta_c(\alpha) < \theta \leq \theta_{-}$.

Following a symmetric reasoning for negative impulses, that is, $u(\cdot) = -\alpha\,\delta(\cdot)$, the phase response curve is given by
\begin{align*}
	\PRC^0(\theta;\alpha\,\delta(\cdot)) &= 
	\begin{cases}
		\theta_{+} + \freqslow^0 \, \psi_{-}(z_{+},z^{\perorb}(\theta),0) - \theta, &\text{if (C4)}, \\
		0, & \text{o/w},
	\end{cases}
\end{align*}
where (C4) stands for $\theta_c(\alpha) < \theta \leq \theta_{+}$ and $z_{c}(\alpha) = \{ z\in\Rbb :  \critfun_{+}(z) - \alpha = \critfun_{r}(z) \}$.

\subsubsection{Square pulse of finite duration}

A square pulse of finite duration $u(\cdot) = \bar{u}\,\left[\ech(\cdot) - \ech(\cdot - \Delta)\right]$ induces a behavior in the fast-slow dynamics \eqref{eq:relax_fast_syst layer}--\eqref{eq:relax_slow_syst reduced} that is less trivial. 

The phase response curve is given by
\begin{align}
	\PRC^0(\theta;u(\cdot)) & = \Theta^0(x_\Delta(\theta),z_\Delta(\theta)) - (\theta + \freqslow^0 \, \Deltaslow^0)
\end{align}
where $(x_\Delta(\theta),z_\Delta(\theta))$ is the state at time $\Deltaslow^0$ for the reduced dynamics starting from $(x^{\perorb}(\theta),z^{\perorb}(\theta))$ where $\Deltaslow^0$ is the pulse duration in the slow time scale and in the singular limit.
It is thus necessary to compute the  state $(x_\Delta,z_\Delta)$ of the trajectory at the end of the pulse in order to compute the reset phase associated with its initial condition. 

\begin{figure*}
	\centering
	\includegraphics[scale=0.8]{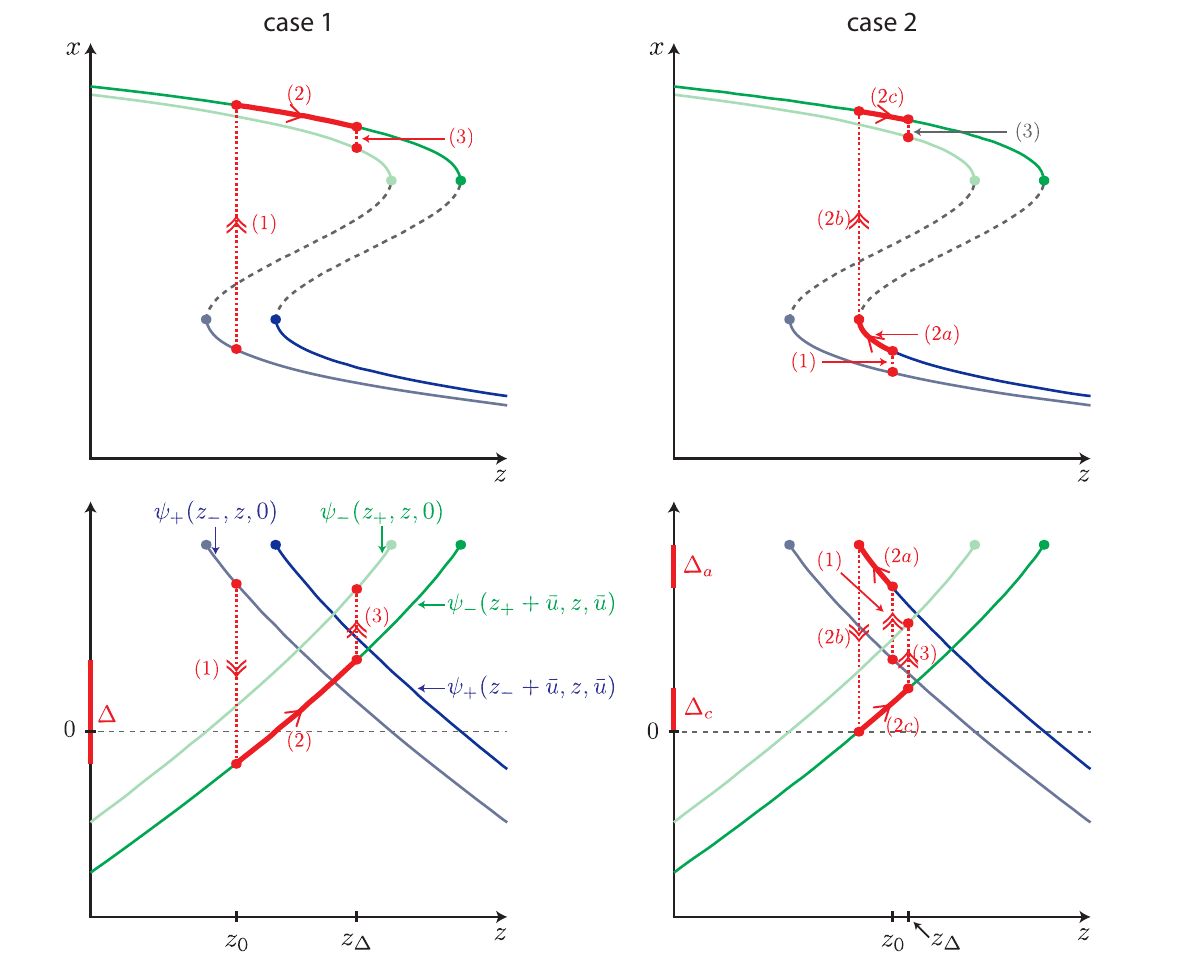}
	\caption{
		Effect of positive square pulses of finite duration in the fast-slow dynamics \eqref{eq:relax_fast_syst layer}--\eqref{eq:relax_slow_syst reduced}. 
		The state $(x_{\Delta},z_{\Delta})$ of the trajectory starting from initial condition $(x_0,z_0)$ (under a pulse of duration $\Delta$) is graphically determined using functions $\psi_{\bullet}$ in order to predict the phase response associated with this pulse. 
		The effect of a positive pulse is to shift temporally the critical manifold along the $z$-axis to the right. The singular trajectory starting from $(x_0,z_0)$ evolves as follows: 
		(1)~jumps instantaneously on the shifted critical manifold, 
		then (2)~evolves around the shifted hysteresis (for a duration $\Delta = \Delta_a + \Delta_c$), and
		finally (3)~jumps back to the initial critical manifold. 
		The main difference between case 1 and case 2 is that during step (1) the trajectory converges to the opposite branch (with respect to the initial point) of the shifted critical manifold (in case 1) or to the same branch  (with respect to the initial point) of the shifted critical manifold (in case 2).
		}
	\label{fig:relax_oscillator_pulse}
\end{figure*}

In the following, we describe how we can compute the  state $(x_\Delta,z_\Delta)$  using only the information contained in the functions $\psi_{-}(z_{+}+\bar{u},z,\bar{u})$ and $\psi_{+}(z_{-}+\bar{u},z,\bar{u})$ (see~\myfigurename~\ref{fig:relax_oscillator_pulse}). 

Starting from the initial condition $(x_0,z_0)$ on the critical manifold, the trajectory evolves as follows (see \myfigurename~\ref{fig:relax_oscillator_pulse}).
\begin{enumerate}  \renewcommand{\labelenumi}{(\arabic{enumi})}
	\item 
	Under a constant input $\bar{u}$, the critical manifold of the system is shifted along the $z$-axis.
	The singular trajectory jumps thus instantaneously to the branch of the  ``shifted critical manifold'' corresponding to the basin of attraction to which the initial state belongs.
	\item Then, the trajectory evolves on the ``shifted critical manifold'', sliding slowly on branches and jumping instantaneously when it reaches ``shifted folds''. 
	\item Finally, the trajectory jumps instantaneously back to the critical manifold at the end of the pulse.
\end{enumerate}
Because the slow variable $z$ is one-dimensional, the evolution of a trajectory under constant input $\bar{u}$ on an attractive branches is fully characterized by the functions $\psi_{-}(z_{+}+\bar{u},z,\bar{u})$ and $\psi_{+}(z_{-}+\bar{u},z,\bar{u})$ during the flowing time. 
The total flowing time has to be equal to the duration $\Deltaslow$. 

In \myfigurename~\ref{fig:relax_oscillator_pulse}, we differentiate between two cases. 
In case~1, the initial condition on the \emph{lower} branch of the critical manifold jumps directly to the \emph{upper} branch of the shifted critical manifold.
In case~2, the initial condition on the \emph{lower} branch of the critical manifold jumps on the lower branch of the ``shifted critical manifold''. 
Case 1 produces larger phase shift than case 2.

\begin{remark} 
	The duration $\Delta$ is expressed in the fast time scale, that is, $\Deltafast^\epsilon = \Delta$. In the slow time scale, the duration is given by $\Deltaslow^\epsilon = \epsilon\,\Deltafast^\epsilon$.  We assume the duration of the pulse $\Deltaslow^\epsilon$ (in the slow time scale) do not tend to zero in the singular limit and thus that the duration $\Deltafast^\epsilon$ tends to infinity. This assumption is motivated by the fact that the duration of the pulse is often a fraction of the period. So we may have $\lim_{\epsilon\rightarrow0}\Deltafast^\epsilon = +\infty$ and $\lim_{\epsilon\rightarrow0}\perfast^\epsilon = +\infty$, and a finite ratio $\lim_{\epsilon\rightarrow0}\Deltafast^\epsilon/\perfast^\epsilon = C$ (with $C\neq0$ and $C\neq\infty$).
\end{remark}

\section{Application to a neural oscillator model} \label{sec:application}

We illustrate our geometric approach on a simple neural oscillator model developed by FitzHugh~\cite{Fitzhugh:1961il} and Nagumo~\cite{Nagumo:1962iz}. This model is a popular two-dimensional simplification of the Hodgkin-Huxley model of spike generation
\begin{subequations}
	\begin{align*}
		\dot{v} & = v - v^3/3 - w + u  \\
		\tau\,\dot{w} & = a - b\,w  + v 
	\end{align*}
\end{subequations}
where $v$ is the voltage variable, $w$ is the recovery variable, and $\epsilon \eqdef 1/\tau$ is a small parameter.

\subsection{Phase response curves for impulses}

\begin{figure}
	\centering
	\includegraphics[scale=0.8]{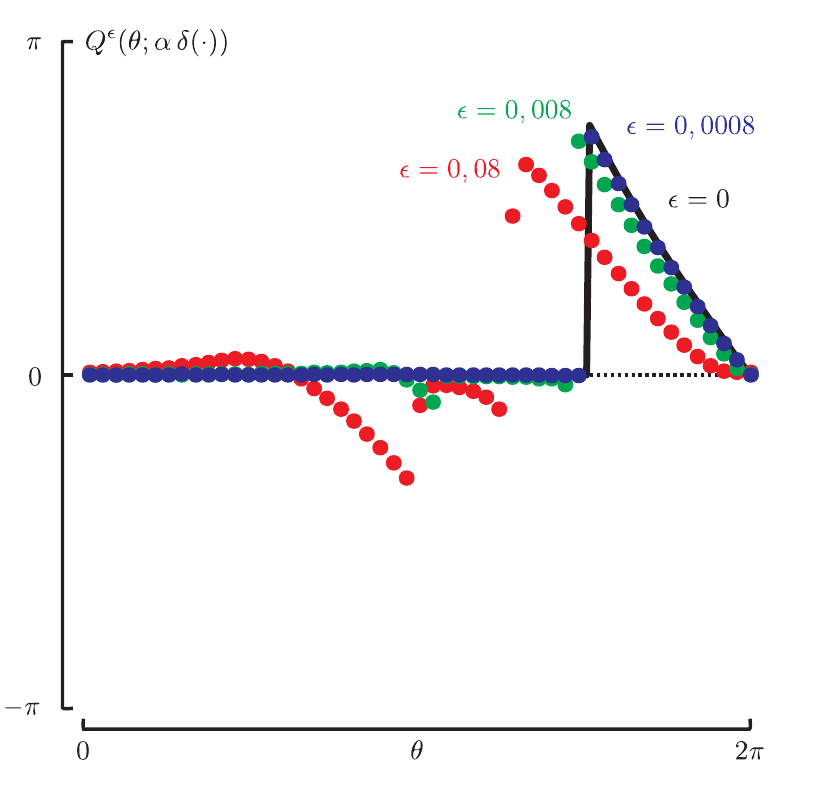}	
	\caption{%
		Phase response curves for excitatory impulses: singular geometric prediction (solid line) and numerical simulations (dots).	
		(Parameter values: $a = 0.7$, $b = 0.8$, $I = 1$, $|\alpha| = 1.5$.)
		}
	\label{fig:FHN-impulse}
\end{figure}

\myfigurename~\ref{fig:FHN-impulse} illustrates the (finite) phase response curve of the FitzHugh-Nagumo model for excitatory impulses $u(\cdot) = \alpha\,\delta(\cdot)$, with $\alpha>0$. The solid line is the geometric prediction computed in the singular limit. Dots represent the phase response computed through numerical simulations of trajectories of the model for different values of the parameter~$\epsilon$. 

The singular phase response curve is equal to zero except in one region of the periodic orbit which corresponds to the region right before the initiation of the upper part of the periodic orbit for an excitatory impulse. In this region, an impulse advances the initiation of the upper part of the periodic orbit. The phase advance decreases monotonically to zero until the phase corresponding to the lower fold.

For small values of $\epsilon$, the geometric prediction matches very well the numerical phase response curves. For larger values of $\epsilon$, the prediction still matches (qualitatively) the larger phase shifts arising before the lower fold but do not capture the small phase shifts arising before the upper fold.

\subsection{Phase response curves for square pulses of finite duration}

\begin{figure}
	\centering
	\includegraphics[scale=0.8]{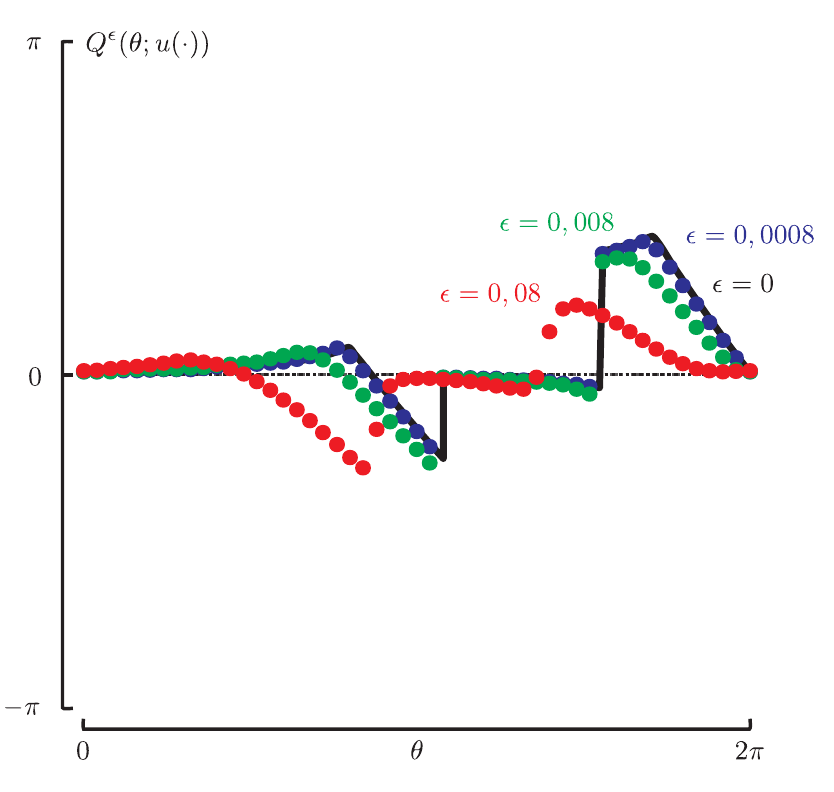}	
	\caption{
		Phase response curves for excitatory pulses of finite duration: singular geometric prediction (solid line) and numerical simulations (dots).	
		(Parameter values: $a = 0.7$, $b = 0.8$, $I = 1$, $|\bar{u}| = 0.25$, $\Delta = 0.1\,T$.)
		}
	\label{fig:FHN-pulse}
\end{figure}

\myfigurename~\ref{fig:FHN-pulse} illustrates the (finite) phase response curve of the FitzHugh-Nagumo model for excitatory square pulses of finite duration. The solid line is the geometric prediction computed in the singular limit. Dots represent the phase response computed through numerical simulations of trajectories of the model for different values of the parameter~$\epsilon$. 

The singular phase response curve is equal to zero except in two regions of the periodic orbit. The first region which exhibits the highest phase shifts corresponds to same region as for the impulse case. The phase shifts in this region follow a piecewise law: the breaking point in the phase shifts corresponds to the separation between initial conditions that continue to evolve on the shifted initial branch  and those that directly jump to the opposite branch.
The second region corresponds to point close to the other fold (see case~1 and case~2 in \myfigurename~\ref{fig:relax_oscillator_pulse}). An excitatory pulse may delay the termination of the upper part.

Once again, for small values of $\epsilon$, the geometric prediction matches very well the actual phase response curves. For larger values of $\epsilon$, the prediction  matches qualitatively both non-zero regions of the phase response curve.

The main difference between the phase response curve for an impulse and for a pulse is that a positive pulse may delay the termination  of the behavior on the upper branch, while a positive impulse may not.

%

\section{Conclusion}

In this paper we overcome the limitation of the infinitesimal phase response curve approach to singularly perturbed oscillators by studying geometrically the asymptotic phase map and the finite phase response curve in the singular limit. 
By persistence results, this analysis provide semi-analytic prediction on the qualitative isochrons structure and on the (finite) phase response curve shape for arbitrary inputs. 
 This result is illustrated on a relaxation oscillator forced by impulses and pulses of finite duration, confirming the goodness of the approach.

Future work will aim at extending this analysis to more complex singularly perturbed oscillators, like bursters \cite{Franci:2013uq}, and to oscillator synchronization studies, linking this result to fast threshold modulation phenomenon \cite{Somers:1993ui}.

%
%

\section*{Acknowledgments}
{
Rodolphe Sepulchre is gratefully acknowledged for insightful comments and suggestions.

This paper presents research results of the Belgian Network DYSCO (Dynamical Systems, Control, and Optimization), funded by the Interuniversity Attraction Poles Programme, initiated by the Belgian State, Science Policy Office. The scientific responsibility rests with its authors. 

P.~Sacr\'e was supported as a Research Fellow with the Belgian Fund for Scientific Research (F.R.S.-FNRS).
}


\bibliographystyle{IEEEtran}
\bibliography{/Users/pierresacre/Dropbox/Research/biblio/biblio_sacre/journal-name-abbrv,/Users/pierresacre/Dropbox/Research/biblio/biblio_sacre/bibfile}

\begin{thebibliography}{10}
\providecommand{\url}[1]{#1}
\csname url@samestyle\endcsname
\providecommand{\newblock}{\relax}
\providecommand{\bibinfo}[2]{#2}
\providecommand{\BIBentrySTDinterwordspacing}{\spaceskip=0pt\relax}
\providecommand{\BIBentryALTinterwordstretchfactor}{4}
\providecommand{\BIBentryALTinterwordspacing}{\spaceskip=\fontdimen2\font plus
\BIBentryALTinterwordstretchfactor\fontdimen3\font minus
  \fontdimen4\font\relax}
\providecommand{\BIBforeignlanguage}[2]{{%
\expandafter\ifx\csname l@#1\endcsname\relax
\typeout{** WARNING: IEEEtran.bst: No hyphenation pattern has been}%
\typeout{** loaded for the language `#1'. Using the pattern for}%
\typeout{** the default language instead.}%
\else
\language=\csname l@#1\endcsname
\fi
#2}}
\providecommand{\BIBdecl}{\relax}
\BIBdecl

\bibitem{Winfree:1980ue}
A.~T. Winfree, \emph{The Geometry of Biological Time}, 1st~ed., ser.
  Biomathematics.\hskip 1em plus 0.5em minus 0.4em\relax New York, NY:
  Springer-Verlag, 1980, vol.~8.

\bibitem{Glass:1988ub}
L.~Glass and M.~C. Mackey, \emph{From Clocks to Chaos: the Rhythms of
  Life}.\hskip 1em plus 0.5em minus 0.4em\relax Princeton, NJ: Princeton
  University Press, 1988.

\bibitem{Efimov:2009fr}
D.~V. Efimov, P.~Sacr{\'e}, and R.~Sepulchre, ``Controlling the phase of an
  oscillator: a phase response curve approach,'' in \emph{Proc. 48th IEEE Conf.
  Decision and Control and 28th Chinese Control Conf.}, Shanghai, China, Dec.
  2009, pp. 7692--7697.

\bibitem{Danzl:2009co}
P.~Danzl, J.~Hespanha, and J.~Moehlis, ``Event-based minimum-time control of
  oscillatory neuron models: phase randomization, maximal spike rate increase,
  and desynchronization,'' \emph{Biol. Cybern.}, vol. 101, no. 5--6, pp.
  387--399, Dec. 2009.

\bibitem{Mauroy:2012vi}
A.~Mauroy, P.~Sacr{\'e}, and R.~Sepulchre, ``Kick synchronization versus
  diffusive synchronization,'' in \emph{Proc. 51st IEEE Conf. Decision and
  Control}, Maui, HI, Dec. 2012, pp. 7171--7183.

\bibitem{Dorfler:2013fk}
F.~Dorfler and F.~Bullo, ``Synchronization in complex oscillator networks: A
  survey,'' Apr. 2013, submitted to Automatica (preprint
  \url{http://motion.me.ucsb.edu/pdf/2013b-db.pdf}).

\bibitem{Sacre:2014aa}
P.~Sacr{\'e} and R.~Sepulchre, ``Sensitivity analysis of oscillator models in
  the space of phase-response curves: Oscillators as open systems,''
  \emph{{IEEE} Control Syst. Mag.}, Apr. 2014, {To appear. Available at
  \href{http://arxiv.org/abs/1206.4144}{arXiv:1206.4144} [math.DS]}.

\bibitem{Izhikevich:2000hb}
E.~M. Izhikevich, ``Phase equations for relaxation oscillators,'' \emph{{SIAM}
  J. Appl. Math.}, vol.~60, no.~5, pp. 1789--1804, 2000.

\bibitem{Krupa:2001ez}
M.~Krupa and P.~Szmolyan, ``Relaxation oscillation and canard explosion,''
  \emph{J. Differential Equations}, vol. 174, no.~2, pp. 312--368, Aug. 2001.

\bibitem{Franci:2013uq}
A.~Franci, G.~Drion, and R.~Sepulchre, ``Modeling the modulation of neuronal
  bursting: a singularity theory approach,'' May 2013, {Available at
  \href{http://arxiv.org/abs/1305.7364}{arXiv:1305.7364} [math.DS]}.

\bibitem{Sacre:2013ys}
P.~Sacr{\'e}, ``Systems analysis of oscillator models in the space of phase
  response curves,'' Ph.D. dissertation, University of Li{\`e}ge, Belgium, Sep.
  2013.

\bibitem{Fitzhugh:1961il}
R.~FitzHugh, ``Impulses and physiological states in theoretical models of nerve
  membrane,'' \emph{Biophys. J.}, vol.~1, no.~6, pp. 445--466, Jul. 1961.

\bibitem{Nagumo:1962iz}
J.~Nagumo, S.~Arimoto, and S.~Yoshizawa, ``An active pulse transmission line
  simulating nerve axon,'' \emph{Proc. IRE}, vol.~50, no.~10, pp. 2061--2070,
  Oct. 1962.

\bibitem{Somers:1993ui}
D.~Somers and N.~Kopell, ``Rapid synchronization through fast threshold
  modulation,'' \emph{Biol. Cybern.}, vol.~68, no.~5, pp. 393--407, Mar. 1993.

\end{thebibliography}

\end{document}